\documentclass[12pt]{amsart}
\usepackage{fullpage,url,amssymb,enumerate,colonequals}

\usepackage{mathrsfs} % for \mathscr (script letters)
\usepackage{MnSymbol}
\usepackage{extarrows}
\usepackage{lscape}
\usepackage[all,cmtip]{xy}
\usepackage{xcolor}

\usepackage[OT2,T1]{fontenc}
\usepackage{color}

\usepackage[
        colorlinks, citecolor=darkgreen,
        backref,
        pdfauthor={Nuno Freitas}, % add other authors
]{hyperref}
\usepackage{cleveref}
\usepackage{comment}

\newcommand{\F}{\mathbb{F}}

\newcommand{\Q}{\mathbb{Q}}

\newcommand{\Z}{\mathbb{Z}}
\newcommand{\Qbar}{{\overline{\Q}}}

\newcommand{\rhobar}{{\overline{\rho}}}

\newcommand{\fp}{\mathfrak p}
\newcommand{\fq}{\mathfrak q}

\newcommand{\calO}{\mathcal{O}}

\newcommand{\Fp}{\mathfrak{p}}

\DeclareMathOperator{\Frob}{Frob}
\DeclareMathOperator{\Gal}{Gal}

\DeclareMathOperator{\Norm}{Norm}

\newcommand{\vv}{\upsilon}

\numberwithin{equation}{section}

\newtheorem{theorem}[equation]{Theorem}
\newtheorem{lemma}[equation]{Lemma}

\newtheorem{proposition}[equation]{Proposition}

\theoremstyle{definition}

\theoremstyle{remark}
\newtheorem{remark}[equation]{Remark}

\definecolor{darkgreen}{rgb}{0,0.5,0}

\author{Nuno Freitas}
\address{Instituto de Ciencias Matem\'aticas, CSIC, Calle Nicol\'as Cabrera 13--15, 28049 Madrid, Spain}
\email{nuno.freitas@icmat.es}

\author{Filip Najman}
\address{University of Zagreb, Faculty of Science, Department of Mathematics,  Bijeni\v{c}ka cesta 30, 10000 Zagreb, Croatia}
\email{fnajman@math.hr}

\title{Two results on $x^r + y^r = dz^p$}

\thanks{The second author was supported by the QuantiXLie Centre of Excellence, a  project co-financed by the Croatian Government and European Union  through the European Regional Development Fund - the Competitiveness and Cohesion Operational Programme (Grant KK.01.1.1.01.0004) and by the Croatian Science Foundation under the project no. IP-2018-01-1313. The first author is partially supported by the PID2019-107297GB-I00 grant of the MICINN (Spain)}

\begin{document}

\begin{abstract} This note proves two theorems regarding Fermat-type equation $x^r + y^r = dz^p$ where
$r \geq 5$ is a prime. Our main result shows that, for infinitely many  integers~$d$, the previous equation has no non-trivial primitive solutions such that $2 \mid x+y$ or $r \mid x+y$, for a set of exponents $p$ of positive density.  We use the modular method with a symplectic argument to prove this result.
\end{abstract}

\maketitle

\section{Introduction}

%\section{A result on signature $(r,r,p)$ using the symplectic argument}
\begin{comment}
Let $A,B,C \in \Z$ non-zero and pairwise coprime.
The \emph{Generalized Fermat Equation}
\begin{equation}
 Ax^p+By^q=Cz^r, \qquad  p,q,r \in \mathbb{Z}_{\geq 2}, \qquad \frac{1}{r}+\frac{1}{q}+\frac{1}{p}<1,
\label{E:GFE}
\end{equation}
has been the center of attention among Diophantine equations since Wiles'~\cite{wiles} groundbreaking proof of Fermat's Last Theorem. A solution~$(a,b,c)$ of~\eqref{E:GFE} is called {\it primitive} if~$\gcd(a,b,c) = 1$ and {\it non-trivial} if $abc\neq 0$. A theorem of Darmon--Granville~\cite{DG} states that if $A,B,C$ and the prime exponents~$p, q, r$ are fixed then there are only finitely many non-trivial primitive solutions to~\eqref{E:GFE}, but more is conjectured:

{\bf Conjecture.} Fix $A,B,C$ as above.
Over all choices of prime exponents $p,q,r$ satisfying $1/r+1/q+1/p<1$ the
equation \eqref{E:GFE} admits only finitely many solutions $(a,b,c)$ such that
$abc \ne 0$ and $\gcd(a,b,c)=1$. (Here solutions like $2^3 + 1^q = 3^2$ are counted only once.)

To a triple of prime exponents~$p,q,r$ as above we call a {\it signature}
and we say that~\eqref{E:GFE}
is a {\it Fermat-type equation of  signature $(r,q,p)$}.
\end{comment}

We consider Fermat-type equations of the form
\begin{equation}
\label{E:rrp}
 x^r + y^r = dz^p,
\end{equation}
where $r, p > 3$ are primes and $d$ is an odd positive integer, not an $r$-th power and $r \nmid d$. A solution~$(a,b,c)$ to~\eqref{E:rrp} is called {\it primitive} if~$\gcd(a,b,c) = 1$ and {\it non-trivial} if $abc\neq 0$.

We will study the above equation via the modular method. Since $d$ is not a $r$-th power, we have $d \neq 1$; this avoids the trivial solution~$(1,0,1)$ which is a well known obstruction. Indeed, the modular method aims to obtain a contradiction with the existence of a solution, which is usually not possible once a solution does exist. More precisely, for $d = 1$, in the elimination stage of the method, we would have to distinguish the mod~$p$ representation of the Frey curve attached to a non-trivial putative solution~$(a,b,c)$ from the mod~$p$ representation of the Frey curve attached to the solution $(1,0,1)$. Currently, with the exception of very few special cases, there are no techniques available to do this. The hypothesis that $d$ is odd avoids the solution $(1,1,1)$ for $d=2$ but it also plays (together with $r \nmid d$) an essential r\^ole in the proof of Theorem~\ref{T:symplectic} below (see Remark~\ref{rem:odd_d}).

Equation~\eqref{E:rrp} has been the focus of various recent works. It has been completely resolved for $(r,d)=(5,3)$, $(13,3)$ and all~$p \geq 2$ in \cite{bcdf, bcddf}, and for $(r,d)=(7,3)$ in~\cite{bcdfn}; moreover, in~\cite{Mocanu} it is shown
for $d=1$ and many values of~$r < 150$ that it admits no non-trivial primitive solutions $(a,b,c)$ with $c$ even.
Under the extra assumptions that $d$ is not divisible by~$p$-th powers and is divisible only by primes $q\not\equiv 1 \pmod r$, an older result
due to Kraus implies that, fixed both~$r$ and $p$,
the set of coefficients~$d$ for which~\eqref{E:rrp} has a non-trivial primitive solution is finite (see~\cite[Th\'eor\`eme~1]{KrausComp}).

Given that we will focus on asymptotic results, i.e. for large enough~$p$, we can also assume that $p$ is large compared to~$r$ and $p \nmid d$. We will prove the following theorem.

\begin{theorem}\label{T:symplectic}
Let $r$ and~$d$ be as above. For a set of prime exponents $p$ of positive density, the equation~\eqref{E:rrp} has no non-trivial primitive solutions~$(a,b,c)$
such that $2 \mid a+b$ or $r \mid a+b$.
\end{theorem}

This result has no constraints on the value of~$r$ and~$d$ beyond those introduced in the first paragraph. Instead we have $2$-adic or $r$-adic restrictions on the solutions, that  naturally occur in the study of Fermat-type equations.
In particular, the condition $r \mid a+b$ is equivalent to~$r \mid c$ and this is analogous to what Darmon calls a {\it first case} solution to the Fermat-type equations  $x^p + y^p =  z^r$ (see \cite[Definition 3.6]{Darmon} and results in \S 3 of {\it loc. cit.}).
A result of similar flavor to ours regarding $x^p + y^p =  z^r$ is given in~\cite{bcdf0}, more precisely, it is shown this equation has no non-trivial primitive solutions with $r \mid ab$ and $2 \nmid ab$. We remark that although in {\it loc. cit.} the proof uses Darmon's higher dimensional abelian varieties instead of Frey elliptic curves, the condition $r \mid ab$ plays a similar role to ours $r \mid a+b$ in that it forces the Frey varieties to have multiplicative reduction at~$r$.

The proof of Theorem~\ref{T:symplectic} was inspired by~\cite{hk} and relies on combining the modular method with the symplectic argument. We refer the reader to~\cite{freKraus1} for a quick introduction to Diophantine applications of the symplectic argument; for a comprehensive study of symplectic criteria and further Diophantine applications we refer to~\cite{freKraus2}.

We will also prove the following theorem about~\eqref{E:rrp} where the $2$-adic and $r$-adic  conditions are replaced by a $p$-adic condition.

\begin{theorem}\label{T:main2} Let $r$ and~$d$ be as above. Then
for  prime exponents $p \equiv \pm 1 \pmod{r}$ large enough
the equation~\eqref{E:rrp} has no non-trivial primitive solutions~$(a,b,c)$
such that $p \mid a+b$.
\end{theorem}

%{\bf Acknowledgments.} We thank Alain Kraus, Luis Dieulefait and Imin Chen and the anonymous referee for their comments. Most of the work on this paper was during the visit of the first author to the University of Zagreb, which was funded by the QuantiXLie Centre of Excellence.

\section{The Frey curve}

We first recall the construction of a Frey curve from~\cite{fre}.
Let $r\geq 5$ be a prime and $\zeta=\zeta_r$ a fixed primitive $r$-th root of unity. Let $K=\Q(\zeta_r)^+$ be the maximal real subfield of $\Q(\zeta_r)$.
Define the polynomials
$$f_{1}:=x^2+(\zeta^{2}+\zeta^{-2})xy+y^2 \quad \text{and} \quad f_{2}:=x^2+(\zeta+\zeta^{-1})xy+y^2$$
and the constants
\[
\alpha = \zeta(1-\zeta)(1-\zeta^{-3}), \quad
 \beta = (1-\zeta)(1-\zeta^{-1}), \quad
 \gamma = -(1-\zeta^2)(1-\zeta^{-2}).
\]
We set also
\[
A=A_{a,b} = \alpha (a+b)^2, \quad B=B_{a,b} = \beta f_{1}(a,b), \quad C=C_{a,b} = \gamma f_{2}(a,b) \]
and define the Frey curve
\begin{equation}\label{eq:FoverK}
F_{a,b} : Y^2 = X(X-A_{(a,b)})(X + B_{(a,b)}).
\end{equation}

By construction we have $A+B=C$ and
the following standard invariants
\begin{eqnarray*}
c_4(F_{a,b}) & = & 2^4(A_{a,b}^2 + A_{a,b}B_{a,b} + B_{a,b}^2), \\
c_6(F_{a,b}) & = & 2^5(2A_{a,b}^3 + 3A_{a,b}^2B_{a,b} - 3A_{a,b}B_{a,b}^2 - 2B_{a,b}^3), \\
\Delta(F_{a,b}) & = & 2^4\left(A_{a,b}B_{a,b}C_{a,b}\right)^2.
\end{eqnarray*}
For $n \in \Z_{>0}$ and $x \in K$ denote by~$\mbox{Rad}_{n}(x)$ the product of the primes in~$K$ not dividing~$n$.

Also let~$\fp_r$ denote the unique prime in~$K$ above~$r$.
\begin{proposition}
\label{prop:conductor}
Let $(a,b,c)$ a primitive solution to~\eqref{E:rrp} such that $\gcd(a,b) =1$.
Then the conductor~$N_F$ of the curve  $F = F_{a,b}$ is of the form
$$ 2^s \fp_r^t \mathfrak{c} \mbox{Rad}_{2r}(a+b)$$
where $\mathfrak{c}$ is a squarefree product of primes in $K$ dividing $c$.

Moreover, $s=1$ if $2 \mid a+b$ and $t=1$ if $r \mid a+b$.
\end{proposition}
\begin{proof}
Since $\alpha, \beta, \gamma$ are of the form $\pm \zeta^s (1 - \zeta^t) (1 - \zeta^u)$, where neither $t$ nor $u$ are $\equiv 0 \pmod{r}$, this means that the only prime dividing $\alpha \beta \gamma$ is $\fp_r$ and $\upsilon_{\fp_r}(\alpha \beta \gamma) = 3$.

Recall the factorisation
\begin{equation}\label{E:factorisation}
 a^r + b^r = (a+b)\Phi_r(a,b) = d c^p
\end{equation}
where $\Phi_r(a,b)$ is the $r$-th cyclotomic polynomial; we have $\gcd(a+b,\Phi_r(a,b))=1$ or $r$ since $a,b$ are coprime. Moreover, the polynomials $f_1(a,b)$ and~$f_2(a,b)$ in the definitions of $B$ and~$C$ are factors
of $\Phi_r(a,b)$ over~$K$; see~\cite[\S 2]{fre} for proofs of these elementary properties (in particular, Corollary~2.2 and the discussion after Corollary~2.6 in {\it loc. cit.}).

Let $\fq \nmid 2r$ be a prime in~$K$. Then $\upsilon_\fq(\Delta(F))=2(\upsilon_\fq(A)+\upsilon_\fq(B) +\upsilon_\fq(C)).$

If $\fq\nmid a^r +  b^r$ we have $\upsilon_\fq(\Delta(F))=0$, so $\fq$ is a prime of good reduction.

Assume $\fq \mid a^r+b^r$. Then $\fq\mid A$ or $\fq\mid \Phi_r(a,b)$ and not both simultaneously.

If $\fq\mid \Phi_r(a,b)$ and $\fq\nmid BC$, then $\fq$ is a prime of good reduction.

Otherwise, either $\fq\mid A$ and $\fq\nmid BC$ or $\fq\nmid A$ and $\fq\mid BC$.
In both cases, we get $v_\fq(\Delta(F))>0$ and $v_\fq(c_4(F))=0$, so $\fq$ is a prime of multiplicative reduction.

We now consider the prime~$\fp_r$ above~$r$. Assume $r \mid a+b$. We have
$$\upsilon_{\fp_r}(c_4) =  4, \quad \upsilon_{\fp_r}(c_6) = 6, \quad \upsilon_{\fp_r}(\Delta) = 10 + 4\upsilon_{\fp_r}(a+b).$$

Since $\upsilon_{\fp_r}(a+b) \geq (r-1)/2$ the equation is non-minimal and after a coordinate change we have $\upsilon_{\fp_r}(c_4) = 0$ and $\upsilon_{\fp_r}(\Delta) > 0$. Thus $F$ has bad multiplicative reduction, i.e. $\upsilon_{\fp_r}(N_F) = 1$.

Finally we consider a prime~$\fq_2 \mid 2$. As 2 is unramified in~$K$ we will use~\cite[Table IV]{papado} to read the conductor at~$\fq_2$ in terms of the $\fq_2$-adic valuations of the standard invariants~$(c_4(F),c_6(F),\Delta(F))$.

Assume $2 \mid a+b$ so $2 \nmid \Phi(a,b)$. Since $p > 3$, the shape of equation~\eqref{E:rrp} implies $8\mid a+b$. Moreover,
$$\upsilon_{\fq_2}(c_4) = 4, \quad \upsilon_{\fq_2}(c_6) = 6, \quad \upsilon_{\fq_2}(\Delta) = 4 + 4\upsilon_{\fq_2}(a+b).$$
The corresponding entries in \cite[Table IV]{papado}
give us Tate case 7 with $\vv_{\fq_2}(N_F)= 4$
or the equation for $F$ is non-minimal as $\upsilon_{\fq_2}(a+b)\geq 3$.
We will apply \cite[Prop.~4]{papado} to show that we are in the non-minimal case; we state it here for convenience.

\begin{proposition}\label{prop:pap4}
Let $K/\Q_2$ be a finite extension with ring of integers~$\calO$ and uniformizer~$\pi$. Let
$W/K$ be an elliptic curve given by a Weierstrass model~$(W)$ with standard invariants $a_i$ and~$b_i$. Assume that $(W)$ is in a Tate case $\geq 7$.

(a) There is $r \in \calO$ such that
\[
 b_8 + 3rb_6 + 3r^2b_4 + r^3b_2 + 3r^4 \equiv 0 \pmod{\pi^5}
\]
(b) Choose any $r$ satisfying (a). Then we are in a Tate case $\geq 8$ if and only if there is $s \in \calO$ satisfying
\[
 a_2 + 3r - s a_1 - s^2 \equiv 0 \pmod{\pi^2}
\]

\end{proposition}

For our curve~$F$, we have $a_1 = a_3 = a_6 =0$, $a_2 = B - A$ and $a_4 = -AB$ and $b_8 = -(AB)^2$. Since $\vv_{\fq_2}(A) \geq 3$, taking $r=0$ the condition in part (a) of Proposition~\ref{prop:pap4} is satisfied and
the condition of part (b) becomes
$B \equiv s^2 \pmod{\fq_2^2}$, that is $B$ is a square mod~$\fq_2^2$.

Let $z=-\zeta_r-\zeta_r^{-1}$. Using $a\equiv -b \pmod{\fq_2^2}$ we get
\begin{align*}
  B&=(2+z)(a^2+(z^2-2)ab+b^2)=(2+z)((a+b)^2-(z^2-4)ab) \\
   &\equiv (z+2)(z^2-4)ab \equiv -b^2(z+2)^2(z-2) \pmod{\fq_2^2}.
\end{align*}
Moreover,
$$2-z =\zeta_r + \zeta_r^{-1}+2 = (\zeta_r^{\frac{r+1}{2}} +\zeta_r^{-\frac{r+1}{2}})^2$$

with $\zeta_r^{\frac{r+1}{2}} +\zeta_r^{-\frac{r+1}{2}} \in \calO_K$ because it is fixed by complex conjugation; thus $B$ is a square mod~$\fq_2^2$ so part (b) is satisfied and we are in the non-minimal case.  After a change of variables we get
\begin{equation}
\upsilon_{\fq_2}(c_4)= \upsilon_{\fq_2}(c_6)=0, \quad \upsilon_{\fq_2}(\Delta(F))=-8+4\upsilon_{\fq_2}(a+b)> 0,
\label{eq:val}
\end{equation}
from which we see that $F$ has multiplicative reduction, i.e $\upsilon_{\fq_2}(N_F)= 1$.
%The last statement follows immediately because the additive primes divide $2r$.
\end{proof}

Recall from~\cite[\S 2]{fre} that $\gcd(a+b,\Phi_r(a,b))=1 $ or $r$ and $\Phi_r(a,b)$ is divisible only by~$r$ and primes $q \equiv 1 \pmod{r}$. Write $d=d_0 d_1$ where a prime $q \mid d_1$ if and only if $q \equiv 1 \pmod{r}$.
\begin{proposition}
For large enough~$p$ the following holds. For all primitive solutions $(a,b,c)$ of~\eqref{E:rrp} with $\gcd(a,b)=1$
the representation~$\bar{\rho}_{F,p}$ is irreducible and modular of weight $2$. Moreover, it's Serre level is~$2^s \fp_r^t \mbox{Rad}_{2r}(d_0 d_1')$ where $\fq \mid d_1'$ if and only if $\fq \mid d_1$ and~$\fq \mid \mathfrak{c}$ where $s$,~$t$ and $\mathfrak{c}$ are as in Proposition~\ref{prop:conductor}.
\label{artinc2}
\end{proposition}
\begin{proof}
From \cite[Theorems 4.3 and 4.4]{fre}, it follows that there exists a constant $C_r$, depending only on~$r$, such that for all $p>C_r$ the representation $\bar{\rho}_{F,p}$ is modular of weight 2 and absolutely irreducible.

Let $\fq \nmid 2r$ be a prime dividing $a+b$. From \Cref{prop:conductor}, we know~$\fq$  is a prime of bad multiplicative reduction.

From~\eqref{E:factorisation} we have $a+b = 2^{sp} \cdot r^k \cdot d_0 \cdot c_0^p$  with $s, k \geq 0$  and~$c_0 \mid c$.
Thus
\begin{equation}\label{E:deltamodp}
 \upsilon_{\fq}(\Delta(F))=2\upsilon_{\fq}(A)=4\upsilon_{\fq}(a+b) =
4\upsilon_{\fq}(d_0) + 4p \upsilon_{\fq}(c_0) \equiv 4\upsilon_{\fq}(d_0) \pmod p.
\end{equation}
When $\upsilon_{\fq}(d_0) \neq 0$ then $4\upsilon_{\fq}(d_0) \not\equiv 0 \pmod p$ for large enough~$p$; further enlarging~$p$ we can assume $p \nmid d$.
Hence the primes $\fq \nmid 2r$ dividing $a+b$ that also divide the Serre level of~$\rhobar_{E,p}$ are precisely those dividing~$d_0$.

Finally, the multiplicative primes $\fq \mid \mathfrak{c}$ divide the Serre level if and only if $\fq \mid d_1$. Indeed, recall that $\fp_r$ is the only prime that can divide both $B$ and~$C$, and that they are both coprime
to $a+b$ as well (and hence to~$d_0$) so both $B$ and $C$ are $p$-th powers times a divisor in~$\calO_K$ of~$d_1 r$. We have
\begin{equation} \label{toref}
\upsilon_{\fq}(\Delta(F)) \equiv 2\upsilon_{\fq}(\Delta(d_1)) \pmod{p},
\end{equation}
hence, for large~$p$, a prime $\fq \mid \mathfrak{c}$ satisfies $\upsilon_{\fq}(\Delta(F)) \not\equiv 0 \pmod p$
if and only if $\fq \mid d_1$.
\end{proof}

\section{Proof of Theorem~\ref{T:symplectic}}

We will use the following lemma several times, so we state it here before proceeding further.

\begin{lemma}[{\cite[Lemme 1.6]{hk}}] \label{lem1.6}
Let $E$ and $E'$ be two elliptic curves over $\Q$ and $p, \ell_1, \ell_2$ three distinct primes with $p>2$. Suppose $E[p]$ and $E'[p]$ are isomorphic. Suppose $E$ and $E'$ have multiplicative reduction at $\ell_i$ and $p$ does not divide $v_{\ell_i}(\Delta (E))$, which implies that $p$ does not divide $v_{\ell_i}(\Delta (E'))$. Then the reduction mod~$p$ of $v_{\ell_1}(\Delta(E))v_{\ell_2}(\Delta(E))v_{\ell_1}(\Delta(E'))v_{\ell_2}(\Delta(E'))$
is a square in $\F_p$.
\end{lemma}

Note that from a non-trivial solution~$(a,b,c)$ to equation~\eqref{E:rrp} with $m=\gcd(a,b) > 1$ and $\gcd(a,b,c)=1$,
we obtain a non-trivial primitive solution
to $x^r + y^r = d' z^p$ with $\gcd(a,b)=1$ where $d' = \frac{d}{m^r}$. Since $d$ is not an $r$-power, the same is true for $d'$. Clearly $d'$ is also odd and $r \nmid d'$, hence we are reduced to the case of solutions satisfying $\gcd(a,b) = 1$. The argument below is inspired by~\cite[Theorem 2.1]{hk} and it will show by contradiction that these latter solutions do not exist for $p$ in a set of primes of density $> 0$.

Suppose there is a non-trivial primitive solution~$(a,b,c)$ to~\eqref{E:rrp} with $a,b$ coprime\footnote{In the paper \cite{fre} where our Frey curve~$F$ was introduced,
a solution $(a,b,c)$ is called primitive when $\gcd(a,b)=1$. Here we decided to use the condition $\gcd(a,b,c)=1$ instead because this is more standard in the context of the generalized Fermat equation $Ax^r + By^q = Cz^p$.}.

By Proposition~\ref{artinc2}, for sufficiently large~$p$, we have $\bar{\rho}_{F,p} \simeq \bar{\rho}_{f,\fp}$, where~$f$ is a Hilbert newform of level $2^s \cdot \fp_r^t \cdot \mbox{Rad}_{2r}(d_0 d_1')$, where $\fp$ is a prime in the field of coefficients $K_f$ of~$f$. Since the field of coefficients of~$F$ is $\Q$, by enlarging $p$ if needed, we can assume that $K_f = \Q$. Moreover, since $d \neq 1$, there is at least one prime $\fq$ in~$K$
dividing $\mbox{Rad}_{2r}(d_0 d_1')$, which is a Steinberg prime of~$f$. Therefore, by the Eichler-Shimura correspondence for~$f$, there is an isogeny class of elliptic curves defined over~$K$ corresponding to~$f$. Let~$E$ denote an elliptic curve in that isogeny class. Since $F$ has full $2$-torsion over~$K$, by further enlarging~$p$ if necessary, we can also assume that~$E$ has full $2$-torsion; see \cite[\S 2.4 and \S 4]{FS} for the previous claims.

From the above, in particular, we have an isomorphism of $G_K$-modules $\phi:F[p]\rightarrow E[p]$.

(1) Suppose $2 \mid a+b$. Let $\fq_2 \mid 2$ and $\fq \mid d$ be primes in~$K$.

Recall that $d=d_0 d_1$ where a prime $q \mid d_1$ if and only if $q \equiv 1 \pmod{r}$.

(1a) Assume that $\fq \mid d_0$. Then the curve~$F$ has multiplicative reduction at $\fq_2$ and~$\fq$ by Proposition~\ref{prop:conductor}; we will apply \Cref{lem1.6}.

Note that $a+b = 2^{sp} \cdot r^k \cdot d_0 \cdot c_0^p$ with $s \geq 1$, therefore
from \eqref{eq:val} and~\eqref{E:deltamodp} we have,
$$\upsilon_{\fq_2}(\Delta(F)) \equiv -8 \pmod p \qquad \text{and} \qquad \upsilon_{\fq}(\Delta(F)) \equiv 4\upsilon_{\fq}(d_0) \not\equiv 0 \pmod p.$$
Therefore, by \Cref{lem1.6} we must have
$$-2\upsilon_{\fq}(d_0)= \upsilon_{\fq_2}(\Delta(E)) \upsilon_{\fq}(\Delta(E)) \text{ in } (\F_p)^*/(\F_{p}^*)^2.$$
Let~$k$ be the number of isogeny classes of elliptic curves over~$K$ with full $2$-torsion and conductor
$2 \cdot \mbox{Rad}_{2r}(d) \cdot \fp_r^t$. For $i = 1,\ldots,k$, we let
$E_i$ be a representative of each isogeny class and set
$$n_i:=-2\upsilon_{\fq}(d_0)\upsilon_{\fq_2}(\Delta(E_i)) \upsilon_{\fq}(\Delta(E_i))$$
which are negative integers as the valuations used in their definition are positive.

Observe that our  result now follows if, for a positive density of primes~$p$, we have that $n_i$ is not a square modulo $p$ for all~$i$. We claim that the set of such primes is non-empty. Then the Dirichlet density theorem guarantees that this set has density $\geq (1/2)^{k}.$

To finish the proof we now prove the claim.
Choose $p\equiv 7 \pmod 8$ such that for all odd primes $q\mid \upsilon_{\fq}(d_0)\prod_{i}n_i$ the condition $\left(\frac{q}{p}\right)=1$ is satisfied. Such primes exist by the Dirichlet density theorem. Let $n_i = - 2^s q_1^{e_1}\ldots q_j^{e_j}$ be the prime factorisation.
Then
$$\left(\frac{n_i}{p}\right)=\left(\frac{-1\cdot 2^s \cdot q_1^{e_1}\ldots q_j^{e_j}}{p}\right)=\left(\frac{-1}{p}\right)\left(\frac{2}{p}\right)^s\cdot 1\cdot \ldots 1=-1,$$
that is, $n_i$ is not a square mod~$p$ for all~$i$, as desired.

(1b) Assume now $\fq \mid d_1$. Thus $\fq \nmid a+b$ and,
after replacing $\fq$ by a conjugate if needed, we can assume
that $\fq \mid \mathfrak{c}$ where $\mathfrak{c}$ is given by Proposition~\ref{prop:conductor}. Thus~$F$ has multiplicative reduction at $\fq_2$ and~$\fq$. From \eqref{toref} we also know that
$\upsilon_{\fq}(\Delta(F)) \equiv 2\upsilon_{\fq}(\Delta(d_1')) \not \equiv 0 \pmod{p}$. The conclusion now follows similar to the previous case where the integers $n_i$ are instead defined by
$n_i:=-\upsilon_{\fq}(d_1')\upsilon_{\fq_2}(\Delta(E_i)) \upsilon_{\fq}(\Delta(E_i))$.

(2) Suppose $r \mid a+b$. Now the curve~$F$ has multiplicative reduction at $\fp_r$ and~$\fq$ by Proposition~\ref{prop:conductor}. Because $\upsilon_r(\Phi_r(a,b)) = 1$ when $r \mid a+b$, we have
$a+b = 2^{sp} \cdot r^{kp-1} \cdot d \cdot c_0^p$. Therefore it follows from the proof of Proposition~\ref{prop:conductor} that a minimal discriminant at~$\fp_r$  for~$F$ satisfies
\[
 \upsilon_{\fp_r}(\Delta) = 10 + 4\upsilon_{\fp_r}(a+b) - 12 = -2 + 4(sp-1)\frac{r-1}{2}  \equiv - 2r \pmod{p}
\]
We now apply \Cref{lem1.6} with primes~$\fp_r$ and~$\fq$. Recall that
$\upsilon_{\fq}(\Delta(F)) \equiv 4\upsilon_\fq(d) \pmod{p}$,
leading to $$\upsilon_{\fq}(\Delta(F)) \upsilon_{\fp_r}(\Delta(F)) = -2 \upsilon_\fq(d) r  \in (\F_p)^*/(\F_{p}^*)^2$$ and
$$n_i = -2 \upsilon_\fq(d) r \cdot \upsilon_{\fq_2}(\Delta(E_i)) \upsilon_{\fp_r}(\Delta(E_i)).$$
Since the $n_i$ are again all negative we complete the proof analogously to case (1).
\qed

\begin{remark}
 \label{rem:odd_d}
 The argument above succeeds due to the negative sign in the definition of the integers $n_i$.
 More precisely, this sign arises due to the
 congruence $\upsilon_{\fq_2}(\Delta(F)) \equiv -8 \pmod p$ in case (1) and $\upsilon_{\fq_r}(\Delta(F)) \equiv -2r \pmod p$ in case (2). We observe these congruences hold only because $d$ is odd and not divisible by~$r$.
\end{remark}

\section{Proof of Theorem~\ref{T:main2}}
As in Theorem~\ref{T:symplectic} we apply the modular method with the Frey curve~$F$.
Moreover, the simplifications at the start of the proof of Theorem~\ref{T:symplectic} also apply.

Namely, let $(a,b,c)$ be a non-trivial solution to~\eqref{E:rrp} such that $p \mid a+b$ and $\gcd(a,b) = 1$.
For primes~$p$ sufficiently large, we obtain  $\bar{\rho}_{F,p} \simeq \bar{\rho}_{E_f,p}$,
where $E_f$ is an elliptic curve with full $2$-torsion associated with a Hilbert
newform~$f$ of level  $2^{s} \cdot \fp_r^t \cdot \mbox{Rad}_{2r}(d_0 d_1')$.

Furthermore, $E_f$ has good reduction  at all primes~$\fp \mid p$ in~$K$ as $p \nmid 2rd$, and $F$ has multiplicative reduction at all~$\fp \mid p$ as $p \mid a+b$ by Proposition~\ref{prop:conductor}.

The assumption $p \equiv \pm 1 \pmod{r}$  implies that
$p$ splits completely in~$K= \Q(\zeta_r)^+$.

In particular, locally at any~$\fp \mid p$ the curves $F$ and~$E_f$ become curves over~$\Q_p$; therefore, by restricting $\bar{\rho}_{F,p} \simeq \bar{\rho}_{E_f,p}$ to decomposition subgroups at~$\fp$, we also have the isomorphism
\begin{equation}
\bar{\rho}_{F,p}|_{\Gal(\Qbar_p /\Q_p)} \simeq \bar{\rho}_{E_f,p} |_{\Gal(\Qbar_p /\Q_p)},
\label{E:frobcong}
\end{equation}
where the representation on the right is the reduction of a crystaline representation and the one on the left is the reduction of a semistable non-crystaline representation. More precisely, the $p$-adic representation $\rho_{F,p}$ when restricted to $D_\fp \simeq \Gal(\Qbar_p /\Q_p)$
is given by the first case in~\cite[\S 2.2.3]{volkov} with $e=1$ whilst
$\rho_{E_f,p}$ restricted to~$D_\fp$ is given by the second case also with $e=1$.
In particular, their semisimplifications are given by
\begin{equation}
\label{E:frobcong2}
(\rho_{E_f,p}|_{D_\fp})^{ss} \simeq \eta_u^{-1} \chi \oplus \eta_u \quad \text{ and } \quad (\rho_{F,p}|_{D_\fp})^{ss} \simeq \eta_{-1}^{a} \chi \oplus \eta_{-1}^{a},
\end{equation}
where the characters $\eta_u$, $\eta_{-1}$ are unramified, $\eta_{-1}$ is quadratic, $a \in \{0,1\}$, $\chi$ is the $p$-adic cyclotomic character
and $a_p(E_f/\Q_p) = u^{-1} p + u$ with $u = \eta_u(\Frob) \in \Z_p^\times$. By putting together the reduction mod~$p$ of the representations in~\eqref{E:frobcong2} with the semisimplification of~\eqref{E:frobcong} we conclude that either
\[
 (\eta_u^{-1}\chi, \eta_u ) \equiv (\eta_{-1}^a\chi, \eta_{-1}^a) \pmod{p} \quad \text{ or } \quad
 (\eta_u^{-1}\chi, \eta_u ) \equiv (\eta_{-1}^a, \eta_{-1}^a\chi) \pmod{p}.
\]
Since $\chi \pmod{p}$ is the mod~$p$ cyclotomic character which is ramified and $\eta_u, \eta_{-1}$ are unramified we must be in the first case. Therefore,
\[
 a_p(E_f/\Q_p) = u^{-1} p + u \equiv u = \eta_u(\Frob) \equiv \eta_{-1}^a(\Frob) = \pm 1 \pmod{p}.
\]
Since $a_\Fp(E_f) = a_p(E_f/\Q_p)$ we have
 \[
  a_\fp(E_f) \equiv \pm 1 \pmod{p} \implies   a_\fp(E_f) = \pm 1 + kp, \quad k \in \Z.
\]
Furthermore, from the Weil bound, we get
\[
 |a_\fp(E_f)| \leq 2 \sqrt{\Norm(\fp)} \implies (\pm 1 + kp)^2 \leq 4 p
\]
since $\Norm(\fp) = p$ because $p$ splits completely in~$K$. The previous inequality does not hold for large~$p$ unless $k=0$, that is $a_\fp(E_f) = \pm 1$.

On the other hand, since~$E_f$ has full $2$-torsion over $\Q$ and reduction modulo a rational prime $q$ of good reduction is injective on $E_f(\Q)_{tors}$, it follows that $a_\ell(E_f)$ is even for all primes~$\ell$ in~$K$ of good reduction; in particular, $a_\fp(E_f)$ is even, a contradiction. \qed

\end{document}